\documentclass{article}
\usepackage{amssymb}

\newcommand{\proof}{{\bf Proof: }}
\newcommand{\example}{{\bf Example: }}
\newcommand{\remark}{{\bf Remark:  }}
\newcommand{\remarks}{{\bf Remarks:  }}
\newcommand{\dimv}{\underline{\dim}\,}
\newcommand{\hb}{\hfill$\square$\medskip}

\newcommand{\SSz}{\scriptstyle}

\newcommand{\supp}{{\rm supp }\,}
\newcommand{\suppb}{{\rm \bf supp }\,}
\newcommand{\ses}[3]
{\mbox{$0 \rightarrow #1 \rightarrow #2 \rightarrow #3 \rightarrow 0$}}
\newtheorem{theorem}{Theorem}[section]
\newtheorem{lemma}[theorem]{Lemma}
\newtheorem{definition}[theorem]{Definition}
\newtheorem{proposition}[theorem]{Proposition}
\newtheorem{corollary}[theorem]{Corollary}

\begin{document}
 
\title{\bf The quantic monoid and degenerate quantized enveloping algebras}
\author{Markus Reineke\\[1ex] BUGH Wuppertal, Gau\ss str. 20, D-42097
Wuppertal, Germany\\ (e-mail: reineke@math.uni-wuppertal.de)}
\date{}
 
\maketitle
 
\begin{abstract}
We study a monoid associated to complex semisimple Lie algebras,
called the quantic monoid. Its monoid ring is shown to be isomorphic to a 
degenerate quantized enveloping algebra. Moreover, we provide normal
forms and a straightening algorithm for this monoid. All these results
are proved by a realization in terms of representations
of quivers, namely as the monoid of generic extensions of a quiver with
automorphism.
\end{abstract}
 
\section{Introduction}\label{intro}

In this paper, we introduce and study the quantic
monoid ${\bf\sf U}$, an object associated to any complex semisimple Lie
algebra $\mathfrak{g}$. Its definition is given in terms of generators and
relations (Definition \ref{defu}), which
can be read off from the Cartan matrix of the Lie algebra $\mathfrak{g}
$.\\[1ex]
We show that the monoid ring ${\bf Q}{\bf\sf U}$ of the quantic monoid can be
viewed as a degenerate quantized enveloping algebra in a natural way. More
precisely, we consider a twisted form ${\cal U}^+_q(\mathfrak{g})$ of the
positive part of the quantized enveloping algebra of $\mathfrak{g}$, which 
can be
specialized at $q=0$. This specialization is isomorphic to ${\bf 
Q}{\bf\sf
U}$ (the Degeneration Theorem \ref{dqg}). Such twisted forms of quantized 
enveloping algebras already appear in \cite{KT} in the type $A$ case, and in 
the Hall algebra approach to quantum groups \cite{Ri}.\\[1ex]
We give several natural normal forms for the elements of ${\bf\sf U}$ in
root-theoretic terms (the Parametrization Theorem \ref{para}), using the 
concept of directed partitions of root systems
introduced in \cite{Re2}. Moreover, we provide a straightening rule for
${\bf\sf U}$ (Proposition \ref{str}), yielding a simple algorithm for 
multiplication of elements in
normal form.\\[1ex]
All these structural results are achieved by realizing the quantic monoid in
terms of quiver representations, namely as the monoid of generic extensions of
a quiver with automorphism (the Realization Theorem \ref{real}). This 
generalizes constructions in \cite{Re1}, which deals with the simply laced 
cases (see also \cite{Re3} 
for a generalization to the Kac-Moody case).\\[2ex]
The step from simply laced to arbitrary root systems shows some surprising 
features:\\[1ex]
The degenerate quantized enveloping algebra is no longer defined by just
specializing suitably twisted quantum Serre relations to $q=0$. Instead, one 
has to
impose additional relations (see Definition \ref{defu}), whose nature is quite 
mysterious from the
algebraic point of view. However, they become completely natural from the 
point of view of quiver representations (see Lemma \ref{natdefrel}).\\[1ex]
It is well known that, in contrast to the case of enveloping algebras, 
there is no embedding of arbitrary quantized enveloping algebras ${\cal 
U}^+(\mathfrak{g})$ into simply laced ones. Our approach shows that this 
becomes true again in the degenerate case. Indeed, Definition \ref{defmassub}, 
in combination with Theorem \ref{real}, shows that an arbitrary quantic monoid 
${\bf\sf U}$ always embeds into a simply laced one. Thus, the same is true for 
degenerate quantized enveloping algebras by Theorem \ref{dqg}.\\[2ex]
Whereas \cite{Re1} is mostly formulated from the point of view of quiver
representations, the results of the present paper are formulated in a purely
root-theoretic language to make them easily accessible. However,
all of them depend entirely on techniques from quiver representation
theory. In particular, section \ref{iso}, which constitutes the technical 
heart of the paper, makes free use of such techniques, for example, the 
structure of Auslander-Reiten quivers \cite{ARS}.\\[2ex]
The paper is organized as follows:\\[1ex]
In section \ref{state}, the quantic monoid is defined, and the structural 
results mentioned above are stated. They are illustrated with a detailed 
discussion of a quantic monoid associated to the root system of type $B_3$ 
(Example \ref{exampleb3}).\\[1ex]
Section \ref{quiv} first recollects several facts on quiver representations 
which 
are used in this paper. Several of them are generalized to the case of a 
quiver together with an automorphism (Lemmas \ref{gks} to \ref{symwse}), which 
is the right framework to formulate the Realization Theorem \ref{real} (for 
similar material, see also \cite{Hu}). Note that the alternative approach to 
non-simply laced root systems via species \cite{Ri_sp} does not apply to the 
present setup, since the geometry of quiver representations is used in an 
essential way, requiring an algebraically closed base field.\\[1ex]
These geometric methods are taken from \cite{Re1}, and are used in section 
\ref{mon} to define the monoid of generic extensions associated to a quiver 
with an automorphism (Definition \ref{defmassub}). Several methods of 
\cite{Re1} are generalized to the present setup (Lemmas \ref{pbwm} to 
\ref{eddec}).\\[1ex]
In section \ref{iso}, the Realization Theorem is proved in the form of Theorem 
\ref{comp}. As noted above, this section makes extensive use of (quiver-) 
representation theoretic techniques. The theorem is first reduced to a 
``straightening rule" (Proposition \ref{straightu}), which is proved by a 
reduction to the rank $2$ case via several intermediate steps (Lemmas 
\ref{red1} to \ref{red5}).\\[1ex]
The efforts of section \ref{iso} are finally rewarded in section \ref{proofs}, 
where all statements of section \ref{state} are easily proved using the 
Realization Theorem and properties of the monoid of generic extensions from 
section \ref{mon}.\\[3ex]
{\bf Acknowledgments:} This paper was written during participation in the TMR 
network ERB RMRX-CT97-0100 ``Algebraic Lie Representations". I would like to 
thank K.~Bongartz, A.~Hubery and P.~Littelmann for helpful remarks.
 
\section{Statement of the results}\label{state}
 
Let $C=(a_{ij})_{i,j\in I}$ be a symmetrizable Cartan matrix of finite type
(see \cite{Lu}), i.e.
\begin{itemize}
\item $a_{ii}=2$ for all $i\in I$ and $a_{ij}\in\{0,-1,-2,\ldots\}$ for all
$i\not=j$ in $I$,
\item there exist $d_i\in{\bf Z}$ for $i\in I$ such that the matrix
$(d_ia_{ij})_{i,j}$ is symmetric,
\item $C$ is positive definite.
\end{itemize}
 
We assume that the $d_i$ are positive and minimal.
Let $(I,\leq)$ be a total ordering of $I$. We will now associate a monoid to
the pair $(C,\leq)$, called the quantic monoid. Its definition might look
quite arbitrary at first sight. But Theorem \ref{dqg} below shows that this
monoid is naturally related to (the quantized enveloping algebra associated
to) the Cartan matric $C$.
 
\begin{definition}[Quantic monoid]\label{defu} Define the quantic monoid ${\bf
\sf U}={\bf \sf
U}(C,\leq)$ as the
monoid with set of
generators $I$ and relations
$$i^pj^qi^rj^s=i^{p+r}j^{q+s}$$
if $i<j$, and $(p,q),(r,s)$ are two consecutive entries in the following list
$L_{ij}$:
\begin{itemize}
\item $L_{ij}=((0,1),(1,0))$ if $a_{ij}=0$,
\item $L_{ij}=((0,1),(1,1),(1,0))$ if $a_{ij}=-1=a_{ji}$,
\item $L_{ij}=((0,1),(1,2),(1,1),(1,0))$ if $a_{ij}=-1, a_{ji}=-2$,
\item $L_{ij}=((0,1),(1,1),(2,1),(1,0))$ if $a_{ij}=-2, a_{ji}=-1$,
\item $L_{ij}=((0,1),(1,3),(1,2),(1,2),(2,3),(1,1),(1,0))$ if $a_{ij}=-1,
a_{ji}=-3$,
\item $L_{ij}=((0,1),(1,1),(3,2),(2,1),(2,1),(3,1),(1,0))$ if $a_{ij}=-3,
a_{ji}=-1$.
\end{itemize}
\end{definition}
 
{\remarks $\;$
\begin{enumerate}
\item In the case of a simply laced Cartan matrix $C$, i.e.~if
$a_{ij}\in\{0,-1\}$ for all $i\not=j$ in $I$, we thus have the following
relations for $i<j$ (see \cite{Re1}):
$$ij=ji\mbox{ if }a_{ij}=0,\;\;\;\left.\begin{array}{ccc}iji&=&iij\\
jij&=&ijj\end{array}\right\}\mbox{ if }a_{ij}=-1.$$
\item The pairs $(p,q)$ in $L_{ij}$ correspond -- via $(p,q)\leftrightarrow
p\alpha_i+q\alpha_j$ -- precisely to the positive roots of the rank $2$ root
system spanned by the simple roots $\alpha_i,\alpha_j$; the only
exception being one root in type $G_2$ {\it which is doubled in $L_{ij}$}.
\item The defining relations of ${\bf \sf U}$ can be rewritten as ``framed
commutation relations" for all $i<j$ in $I$:
\begin{itemize}
\item $[i,j]=0$ if $a_{ij}=0$,
\item $i[i,j]=0$, $[i,j]j=0$ if $a_{ij}=-1$,
\item $i[i,j]=0$, $ij[i,j]j=0$, $[i,j]j^2=0$ if $a_{ij}=-1$, $a_{ji}=-2$,
\item $i[i,j]=0$, $i^2j^2[i,j]j=0$, $ij[i,j]j^2=0$, $ij^2[i,j]j^2=0$,
$[i,j]j^3=0$, $ij[i,j]ij^3=0$ if $a_{ij}=-1$, $a_{ji}=-3$.
\end{itemize}
(The last two cases have obvious dual analogues). The equivalence of these
sets of defining relations can be verified by an elementary calculation.
\item The defining relations of ${\bf \sf U}$ can also be rewritten as
commutation relations
$$[i^pj^q,i^rj^s]=0$$
if $i<j$, and $((p,q),(r,s))$ are two consecutive entries in $L_{ij}$, {\it
together with the relation $ijij^3=ij^2ij^2$ (resp.~$i^3jij=i^2ji^2j$)
in the $G_2$ cases}. Again, the equivalence of the defining relations can be
verified by an elementary calculation. This reformulation is related to a
straightening rule discussed at the end of this section.
\end{enumerate}}

We introduce some basic notation related to quantized enveloping
algebras.\\[1ex]
Let $R=R(C)$ be the root system corresponding to $C$, and let $Q=Q(R)$ be the
root lattice, which we identify with ${\bf Z}I$ via $\alpha_i\leftrightarrow
i$. Similarly, we denote by $R^+$ the set of positive roots, and by $Q^+\simeq
{\bf N}I$ the positive
span of $R^+$ in $Q$.\\[1ex]
Let ${\cal U}^+_v$ be the positive part of the quantized enveloping algebra
(over ${\bf Q}(v)$) associated to $C$ (see \cite{Lu}); it is given by
generators $E_i$ for $i\in I$ and the quantized Serre relations
$$\sum_{p+p'=1-a_{ij}}(-1)^{p'}\left[{{p+p'}\atop
p}\right]_{d_i}E_i^pE_jE_i^{p'}=0\mbox{
for all }i\not=j\mbox{ in }I,$$
where $\left[{{p+p'}\atop p}\right]_{i}$ denotes the usual quantum binomial
coefficients defined via the quantum numbers
$[n]_{i}:=(v_{i}^n-v_{i}^{-n})/(v_{i}-v_{i}^{-1})$, $v_{i}:=v^{d_i}$ (see
\cite{Lu}).\\[1ex]
The algebra ${\cal U}^+_v$ is $Q^+$-graded by defining the degree of $E_i$ as
$\alpha_i\in Q^+$. The degree of a homogeneous element $x\in{\cal U}^+_v$ is
denoted
by $|x|$.\\[1ex]
Using the ordering on $I$, we can consider the (non-symmetric) inner product
$\langle,\rangle$ on $Q$ given by
\begin{equation}\label{eulerform} \langle 
\alpha_i,\alpha_j\rangle=\left\{\begin{array}{ccc} d_ia_{ij}&,&i<j\\
d_i&,&i=j\\
0&,&i>j\end{array}\right. .\end{equation}
Note that the symmetrization of this form is the symmetric form
$(d_ia_{ij})_{ij}$, which is a Cartan datum in the sense of (\cite{Lu},
1.1.1.). The bilinear form  $\langle,\rangle$ allows us to define a variant of
${\cal U}_v^+$ as follows:
 
\begin{definition} Define a new multiplication $*$ on ${\cal U}^+_v$ by
$$x*y=v^{-\langle|x|,|y|\rangle}x\cdot y$$
on homogeneous elements $x,y\in{\cal U}^+_v$.
Set $q=v^2\in{\bf Q}(v)$ and define ${\cal U}^+_q$ as the ${\bf
Q}[[q]]$-subalgebra generated by the $E_i$ for $i\in I$.
\end{definition}
 
\remarks $\;$
\begin{enumerate}
\item In other words, we produce a version of ${\cal U}^+_v$ which can be
specialized at $q=0$ by ``breaking the symmetry" of ${\cal U}^+_v$ and
considering a natural ${\bf Q}[[q]]$-subalgebra.
\item A short calculation using the definition of the twisted multiplication
$*$ shows that
the algebra ${\cal U}^+_q$ fulfills the ``de-symmetrized Serre relations" (see
\cite{Ri}) given in the lemma below. But in general, these relations are no
longer defining. Examples for this can
be seen in the proof of Theorem \ref{dqg} in section \ref{proofs}.
\end{enumerate}
 
\begin{lemma}\label{qserre} The following relations hold in ${\cal U}^+_q$:
\begin{eqnarray*}
\sum_{p+p'=1-a_{ij}}(-1)^{p'}\left\{{{p+p'}\atop
p}\right\}_{i}q_i^{p'(p'-1)/2}E_i^pE_jE_i^{p'}&=&0,\\
\sum_{p+p'=1-a_{ji}}(-1)^{p}\left\{{{p+p'}\atop
p}\right\}_{j}q_j^{p(p-1)/2}E_j^pE_iE_j^{p'}&=&0
\end{eqnarray*}
for all $i<j$ in $I$, where $q_i=q^{d_i}$, and the quantum binomial
coefficient is defined via the quantum numbers
$\{n\}_i=(q_i^n-1)/(q_i-1)$.
\end{lemma}

We can now formulate the first main result of this paper:
 
\begin{theorem}[Degeneration]\label{dqg} The specialization ${\bf 
Q}\otimes_{{\bf 
Q}[[q]]}{\cal U}^+_q$ of ${\cal U}^+_q$
at $q=0$ is isomorphic to
the monoid ring ${\bf Q}{\bf \sf U}$.
\end{theorem}

In other words, the monoid ring ${\bf Q}{\bf \sf U}$ can be viewed as a
degenerate quantized enveloping algebra. In particular, the theorem justifies
the -- at first sight -- complicated defining relations of ${\bf \sf
U}$.\\[2ex]
The second result concerns an explicit realization of ${\bf \sf U}$ in terms
of quiver representations. In fact, this realization is the basis for the
proof
of all other statements in this section. The precise formulation of the
theorem will be postponed to section \ref{iso} (see Theorem \ref{comp}).
 
\begin{theorem}[Realization]\label{real} The quantic monoid ${\bf \sf
U}(C,\leq)$ is isomorphic to the monoid of generic extensions of a quiver with
automorphism associated to $(C,\leq)$.
\end{theorem}
 
Choosing an enumeration $i_1<i_2<\ldots<i_n$ of $I$, we can view any element
$d=\sum_{i\in I}d_i\alpha_i$ in $Q^+$ as an element of ${\bf \sf U}$ via
$$(d)=i_1^{d_{i_1}}\cdot\ldots\cdot i_n^{d_{i_n}}\in{\bf \sf U}.$$
To analyse the behaviour of these natural elements of ${\bf \sf U}$ under
multiplication, we need some additional notation.\\[2ex]
Let ${\bf N}R^+$ be the set of functions on $R^+$ with values in ${\bf N}$.
Define the weight $|{\bf a}|$ of a function ${\bf a}\in{\bf N}R^+$ as $|{\bf
a}|=\sum_{\alpha\in R^+}{\bf a}(\alpha)\alpha\in Q^+$. Using the decomposition
$d=\sum_{\alpha\in R^+}a_\alpha\alpha$ provided by the following lemma, we
associate to any $d\in Q^+$ a function ${\bf a}_d\in{\bf N}R^+$ by ${\bf
a}_d(\alpha)=a_\alpha$. We have $|{\bf a}_d|=d$.
 
\begin{lemma}\label{udec} Given $d\in Q^+$, there exists a unique
de\-com\-po\-sition $$d=\sum_{\alpha\in R^+}a_\alpha\alpha$$ such that
$\langle \alpha,\beta\rangle\geq 0$ if $a_\alpha\not=0\not=a_\beta$.
\end{lemma}
 
Using this notation, we can formulate:
 
\begin{proposition}[Straightening rule]\label{str} Given $d,e\in Q^+$, assume
that
$$\langle \beta,\alpha\rangle\geq 0\mbox{ for all }\alpha,\beta\in R^+\mbox{
such that }{\bf a}_d(\alpha)\not=0\not={\bf a}_e(\beta).$$
Then $$(d)\cdot(e)=(d+e)\mbox{ in } {\bf \sf U}.$$
\end{proposition}
 
Finally, we will construct several parametrizations of the elements of ${\bf
\sf U}$ by using the concept of a directed partition (see \cite{Re2},
\cite{Re4}).
 
\begin{definition}\label{directed} Define a directed partition ${\cal I}_*$ of 
$R^+$ to be a
partition into disjoint subsets $R^+={\cal I}_1\cup\ldots\cup{\cal I}_k$ such
that
\begin{enumerate}
\item $\langle \alpha,\beta\rangle\geq 0$ for $\alpha,\beta\in{\cal I}_s$,
$1\leq s\leq k$,
\item $\langle \alpha,\beta\rangle\geq0\geq\langle \beta,\alpha\rangle$ for
all $\alpha\in{\cal I}_s$, $\beta\in{\cal I}_t$, $1\leq s<t\leq k$.
\end{enumerate}
\end{definition}
 
\begin{lemma}\label{ede} There exists an enumeration
$\alpha^1,\ldots,\alpha^\nu$ of
$R^+$ such that $\langle \alpha^k,\alpha^l\rangle\geq 0$ if $k\leq l$ and
$\langle \alpha^l,\alpha^k\rangle\leq 0$ if $k>l$.
\end{lemma}
 
\remark This lemma shows that directed partitions do exist: given an 
enumeration
as above, $R^+=\{\alpha^1\}\cup\ldots\cup\{\alpha^\nu\}$ is obviously a
directed partition.\\[1ex]
Fix a directed partition ${\cal I}_*$. We associate an element of ${\bf \sf
U}$
to any function ${\bf a}\in{\bf N}R^+$ by
$$({\bf a})=(\sum_{\alpha\in{\cal I}_1}{\bf a}(\alpha)\alpha)\cdot\ldots
\cdot(\sum_{\alpha\in{\cal I}_k}{\bf a}(\alpha)\alpha)\in{\bf \sf U}.$$
 
\begin{theorem}[Parametrization]\label{para} The map $$p_{{\cal I}_*}:{\bf 
N}R^+\rightarrow
{\bf \sf
U},\;\;\; {\bf
a}\mapsto({\bf a})$$
is a bijection.
\end{theorem}
 
\remark In other words, relative to a directed partition, we get a
parametrization of the elements of ${\bf \sf U}$, as well as a normal form for
them. In fact, we will see in section \ref{mon} that Proposition \ref{str}
(or, more
precisely, a special case of it) can be viewed as a straightening rule, which
allows us to straighten an arbitrary word in the alphabet $I$ to the form
provided by the theorem.\\[2ex]
\example\label{exampleb3} We illustrate the results of this section in a 
particular example. 
Let $C$ be the Cartan matrix of type $B_3$ over the index set $I=\{1<2<3\}$. 
The matrix $C$ and the matrix
representing
the non-symmetric form (\ref{eulerform}) are thus given, respectively, by
$$\left(\begin{array}{rrr}2&-1&0\\ -1&2&-2\\ 0&-1&2\end{array}\right)\mbox{ 
and }
\left(\begin{array}{rrr}1&-1&0\\ 0&1&-2\\ 0&0&2
\end{array}\right).$$
The following diagram gives the positive roots, where $(l_1,l_2,l_3)$ denotes 
the root $l_1\alpha_1+l_2\alpha_2+l_3\alpha_3\in R^+$. They are presented in 
the form of a graph (in fact, a ``symmetrized" Auslander-Reiten quiver), such 
that there exists a path from $\alpha$ to $\beta$ if $\langle 
\alpha,\beta\rangle>0$ or $\langle\beta,\alpha\rangle<0$. Thus, reading the 
diagram from the left to the right gives an ordering as in Lemma \ref{ede}:
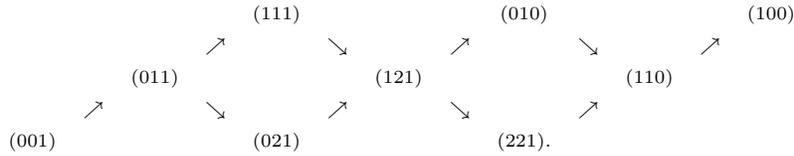
$${\begin{array}{ccccccccccccc}
&&&&\SSz(111)&&&&\SSz(010)&&&&\SSz(100)\\
&&&\SSz\nearrow&&\SSz\searrow&&\SSz\nearrow&&\SSz\searrow&&\SSz\nearrow&\\
&&\SSz(011)&&&&\SSz(121)&&&&\SSz(110)&&\\
&\SSz\nearrow&&\SSz\searrow&&\SSz\nearrow&&\SSz\searrow&&\SSz\nearrow&&&\\
\SSz(001)&&&&\SSz(021)&&&&\SSz(221).&&&&
\end{array}}$$
The quantic monoid ${\bf\sf U}$ has generators $1$, $2$ and $3$,
subjected to the defining relations
$$121=112,\; 212=122,\; 13=31,\; 323=233,\; 23223=22233,\; 2232=2223.$$
Applying Theorem \ref{para} to the trivial directed partition provided by the 
above enumeration of $R^+$, we see that each element of ${\bf\sf U}$ can be 
written as
\begin{equation}\label{normalform}
3^a\cdot(23)^b\cdot(123)^c\cdot(223)^d\cdot(1223)^e\cdot2^f\cdot(11223)^g\cdot(12)^h\cdot 
1^i\end{equation}
for $a,\ldots,i\in{\bf N}$.
The following table gives straightening rules for all the root elements 
$(\alpha)$ for $\alpha\in R^+$. The entry at position $(\alpha),(\beta)$ in 
the table gives a rewriting of the product $(\alpha)\cdot(\beta)$ in ${\bf\sf 
U}$ in the form (\ref{normalform}). All relations are easily verified 
using Proposition \ref{str}.
$${\begin{array}{|c||c|c|c|c|c|c|c|c|} \hline
&\SSz 3&\SSz 23&\SSz 123&\SSz 223&\SSz 1223&\SSz 2&\SSz 11223&\SSz 12\\ 
\hline\hline
\SSz 3&&&&&&&&\\ \hline
\SSz 23&\begin{array}{c}\SSz 3\cdot\\ \SSz 23\end{array}&&&&&&&\\ \hline
\SSz 123&\begin{array}{c}\SSz 3\cdot\\ \SSz 
123\end{array}&\begin{array}{c}\SSz 23\cdot\\ \SSz 123\end{array}&&&&&&\\ 
\hline
\SSz 223&\SSz (23)^2&\begin{array}{c}\SSz 23\cdot\\ \SSz 
223\end{array}&\begin{array}{c}\SSz 123\cdot\\ \SSz 223\end{array}&&&&&\\ 
\hline
\SSz 1223&\begin{array}{c}\SSz 23\cdot\\ \SSz 
123\end{array}&\begin{array}{c}\SSz 123\cdot\\ \SSz 
223\end{array}&\begin{array}{c}\SSz 123\cdot\\
\SSz 1223\end{array}&
\begin{array}{c}\SSz 223\cdot\\ \SSz 1223\end{array}&&&&\\ \hline
\SSz 2&\SSz 23&\SSz 223&\SSz 1223&\begin{array}{c}\SSz 223\cdot\\ \SSz 
2\end{array}&\begin{array}{c}\SSz 1223\cdot\\ \SSz 2\end{array}&&&\\ \hline
\SSz 11223&\SSz (123)^2&\begin{array}{c}\SSz 123\cdot\\ \SSz 
1223\end{array}&\begin{array}{c}\SSz 123\cdot\\ \SSz 11223\end{array}&\SSz 
(1223)^2&
\begin{array}{c}\SSz 1223\cdot\\ \SSz 11223\end{array}&\begin{array}{c}\SSz 
2\cdot\\ \SSz 11223\end{array}&&\\ \hline
\SSz 12&\SSz 123&\SSz 1223&\SSz 11223&\begin{array}{c}\SSz 1223\cdot\\ \SSz 
2\end{array}&\begin{array}{c}\SSz 2\cdot\\ \SSz 11223\end{array}
&\begin{array}{c}\SSz 2\cdot\\ \SSz 12\end{array}&
\begin{array}{c}\SSz 11223\cdot\\ \SSz 12\end{array}&\\ \hline
\SSz 1&\begin{array}{c}\SSz 3\cdot\\ \SSz 1\end{array}&\SSz 
123&\begin{array}{c}\SSz 123\cdot\\ \SSz 1\end{array}&\SSz 1223&\SSz 
11223&\SSz 12&
\begin{array}{c}\SSz 11223\cdot\\ \SSz 1\end{array}&
\begin{array}{c}\SSz 12\cdot\\ \SSz 1\end{array}\\
\hline\end{array}}$$
The subsets
\begin{description}
\item ${\cal I}_1=\{(001),(011),(111)\}$,
\item ${\cal I}_2=\{(021), (121), (010)\}$,
\item ${\cal I}_3=\{(221), (110), (100)\}$
\end{description}
form a directed partition of $R^+$. Theorem \ref{para} implies that any 
element of ${\bf \sf U}$ can be written in the
form
$$1^c2^{b+c}3^{a+b+c}1^{e}2^{2d+2e+f}3^{d+e}1^{2g+h+i}2^{2g+h}3^g$$
for $a,\ldots i\in{\bf N}$. In other words, the set of monomials
$$1^{x_1}2^{x_2}3^{x_3}1^{x_4}2^{x_5}3^{x_6}1^{x_7}2^{x_8}3^{x_9}$$
such that
$$0\leq x_1\leq x_2\leq x_3,\; 0\leq x_4\leq x_6, 2x_6\leq x_4,\; 0\leq x_9, 
2x_9\leq x_8\leq x_7$$
gives a parametrization of the elements of ${\bf\sf U}$.

\section{Quivers with au\-to\-mor\-phisms and their
representations}\label{quiv}
 
Let ${\Gamma}$ be a quiver, i.e.~a finite oriented graph with set of vertices
${\Gamma}_0$.
Let $\gamma$ be an automorphism of ${\Gamma}$, i.e.~a bijection
$\gamma:{\Gamma}_0\rightarrow {\Gamma}_0$ such that for all $i,j\in
{\Gamma}_0$, there is an arrow
from $i$ to $j$ if and only if there is an arrow from $\gamma i$ to $\gamma
j$.\\[1ex]
We always assume ${\Gamma}$ to be of Dynkin type, i.e.~the unoriented
graph
underlying ${\Gamma}$ is assumed to be a disjoint union of Dynkin diagrams of
type
$A$, $D$ and $E$.
A case by case analysis shows that, if ${\Gamma}$ is connected and $\gamma$ is
not
the
identity, then there are precisely the following possibilities for the pair
$({\Gamma},\gamma)$ (see \cite{Lu}):
\begin{itemize}
\item ${\Gamma}$ of type $A_{2n-1}$, $\gamma$ of order $2$ (type $C_n$),
\item ${\Gamma}$ of type $D_{n+1}$, $\gamma$ of order $2$ (type $B_n$),
\item ${\Gamma}$ of type $D_4$, $\gamma$ of order $3$ (type $G_2$),
\item ${\Gamma}$ of type $E_6$, $\gamma$ of order $2$ (type $F_4$).
\end{itemize}
 
We associate to the quiver $\Gamma$ a Cartan matrix
$\widetilde{C}=(\widetilde{a}_{ij})_{i,j\in \Gamma_0}$ by defining 
$-\widetilde{a}_{ij}$
as the number of arrows between $i$ and $j$ (in either direction) for
$i\not=j$. As in the previous section, we denote by $\widetilde{R}^+$ and
$\widetilde{Q}^+$ the corresponding set of positive roots and the positive
part of the root lattice, respectively.\\[1ex]
Moreover, we associate to the pair $(\Gamma,\gamma)$ a Cartan matrix over a
totally
ordered index set as follows (see also (\cite{Lu}, 14.1.1.):\\[1ex]
Let $I$ be the set of $\gamma$-orbits ${\bf i}$ in $\Gamma_0$, and choose a
total ordering on $I$ such that ${\bf i}<{\bf j}$ if there exists an arrow
from a vertex $i\in{\bf i}$ to a vertex $j\in{\bf j}$; such an ordering
exists,
since $\Gamma$, being a Dynkin quiver, has no oriented cycles.
\begin{definition}\label{cartan} For ${\bf i}\not={\bf j}$, define
$-a_{{\bf
i},{\bf j}}$ as the number of arrows between some vertex in ${\bf i}$ and some
vertex in ${\bf j}$ (in either direction), divided by the cardinality of ${\bf
i}$.
\end{definition}
It is then easy to see that $C=(a_{{\bf i}{\bf j}})_{{\bf i},{\bf j}\in
I}$ is a symmetrizable Cartan matrix, with the cardinality $d_{\bf i}$ of the
orbit ${\bf
i}$ as symmetrization index.\\[1ex]
Using this notation, we can identify $Q^+$ with $(\widetilde{Q}^+)^\gamma$,
the $\gamma$-fixed elements in $\widetilde{Q}^+$, via $\alpha_{\bf
i}\leftrightarrow \sum_{i\in{\bf i}}\alpha_i$. This induces an identification
of $R^+$ with the $\gamma$-symmetrizations of elements of $\widetilde{R}^+$.
In the following, we will freely use these identifications; in particular, we
will not distinguish between $\gamma$-fixed elements $d\in
(\widetilde{Q}^+)^\gamma$ and their induced elements in $Q^+$.\\[2ex]
Let $k$ be an algebraically closed field, and denote by $\bmod k{\Gamma}$ the 
category of finite
dimensional
$k$-re\-pre\-sen\-ta\-tions of ${\Gamma}$. This is an abelian $k$-linear
category, since it
is equivalent to the category of finite dimensional representations of the
path algebra $k{\Gamma}$ of ${\Gamma}$ over $k$ (see \cite{ARS}). The 
Auslander-Reiten quiver of $\Gamma$ encodes
the structure of this category: its vertices correspond to the isomorphism 
classes of indecomposable
representations
in $\bmod k\Gamma$, its arrows are given by irreducible maps, and there is an 
additional graph endomorphism $\tau$
corresponding to
the Auslander-Reiten translation (see \cite{ARS} for details).\\[1ex]
For representations $M,N\in \bmod k{\Gamma}$, we denote by $[M,N]$
(resp.~$[M,N]^1$)
the $k$-dimension of ${\rm Hom}_{k{\Gamma}}(M,N)$ (resp.~of ${\rm
Ext}^1_{k{\Gamma}}(M,N)$). We have $[\gamma
M,\gamma
N]=[M,N]$.\\[1ex]
We define a non-symmetric bilinear form $\langle,\rangle$ (the Euler form) on 
${\bf N}\Gamma_0$
by $\langle i,i\rangle=1$, and $\langle i,j\rangle$ equals the number of
arrows from $i$ to $j$ in $\Gamma$ for $i\not=j$. Denoting by $\dimv
M=\sum_{i\in\Gamma_0}\dim_kM_ii\in{\bf N}\Gamma_0$ the dimension vector of a
representation $M$, we have
$$\langle\dimv M,\dimv N\rangle=[M,N]-[M,N]^1.$$
Via the identification of $\widetilde{Q}$ and ${\bf N}\Gamma_0$, the
restriction of $\langle,\rangle$ to $\widetilde{Q}^\gamma\simeq Q$ identifies
with the non-symmetric bilinear form (\ref{eulerform}) on $Q$ introduced in 
the previous section.\\[1ex]
The graph automorphism $\gamma$ of ${\Gamma}$ induces an algebra automorphism
$\gamma$ of $k{\Gamma}$, which is neccessarily of finite order (since $\gamma$
is
so). Given a representation $M\in \bmod k{\Gamma}$, we define a new
representation
$\gamma M\in\bmod k{\Gamma}$ with the same underlying $k$-vector space and the
twisted multiplication $a*m=\gamma(a)m$ for $a\in k{\Gamma}$, $m\in M$. We
define $M$
to be $\gamma$-symmetric if $\gamma M\simeq M$. In this case, $M$ is called
$\gamma$-indecomposable if $M$ has no $\gamma$-symmetric direct summands
except $0$ and $M$ itself.\\[2ex]
For a vertex $i\in \Gamma_0$, denote by $E_i$, $P_i$, $I_i$ the simple,
resp.~indecomposable projective, resp.~indecomposable injective representation
associated to $i$. For an orbit ${\bf i}\in I$, denote by $E_{\bf i}$ the
$\gamma$-symmetric representation $\oplus_{i\in{\bf i}}E_i$; define $P_{{\bf
i}}$ and $I_{{\bf i}}$ similarly. The support $\supp M$ of a representation 
$M\in\bmod k\Gamma$ is the full subquiver
of $\Gamma$ on all vertices $i\in \Gamma_0$ such that $E_i$ appears as a 
composition factor in $M$, or, equivalently,
such
that $P_i$ maps to
$M$. If $M$ is $\gamma$-symmetric, we define its symmetrized support $\suppb 
M$ as the set of all ${\bf i}\in I$ such
that $i\in\supp M$
for some $i\in{\bf i}$.
 
\begin{definition}\label{symmion} Let $U\in\bmod k{\Gamma}$ be an
indecomposable
representation. Let $n\geq 1$ be minimal such that $\gamma^nU\simeq U$; this 
number is called the symmetrization index of $U$. Then
we define the $\gamma$-symmetrization $\widehat{U}$ of $U$ by
$\widehat{U}=\bigoplus_{k=0}^{n-1}\gamma^kU$.
Note that the representation $\widehat{U}$ is obviously $\gamma$-symmetric.
\end{definition}
 
\begin{lemma}\label{gks} $\;$
\begin{enumerate}
\item Any $\gamma$-symmetric representation is a direct sum of
$\gamma$-in\-de\-com\-posables, and this decomposition is unique up to
isomorphisms and permutations.
\item The $\gamma$-indecomposables are of the form $\widehat{U}$ for
indecomposables $U\in\bmod k{\Gamma}$.
\end{enumerate}
\end{lemma}
 
\proof Part a) follows immediately from the definitions and the Krull-Schmidt
theorem for $\bmod k{\Gamma}$. To prove part b), let $M$ be
$\gamma$-indecomposable, and let $U$ be an indecomposable direct summand of
$M$.
Since $\gamma M\simeq M$, all $\gamma^kU$ are again direct summands of $M$.
Thus, $\widehat{U}$ is a direct summand of $M$ which is symmetric, and we
conclude that $M\simeq\widehat{U}$ by $\gamma$-indecomposability of $M$. \hb
 
In particular, for each $\gamma$-orbit ${\bf i}$, we have $E_{\bf 
i}=\widehat{E_i}$ for each vertex $i\in{\bf i}$,
and similarly for $P_{\bf i}$ and $I_{\bf i}$.\\[2ex]
By Gabriel's theorem (see e.g.~\cite{ARS}, VIII), the isomorphism classes of 
indecomposables in $\bmod k\Gamma$ correspond bijectively to the positive 
roots $\widetilde{R}^+$. Moreover, there
exists a partial ordering $\preceq$ on the isomorphism classes of
indecomposables in  $\bmod k{\Gamma}$ such that $[U,V]\not=0$ or
$[V,U]^1\not=0$
implies $U\preceq V$. This ordering can be defined by setting $U\preceq V$ if
there exists a chain of non-zero maps $U=U_0\rightarrow
U_1\rightarrow\ldots\rightarrow U_n=V$, where all $U_k$ are
indecomposable.\\[1ex]
These results extend to the $\gamma$-symmetrized case. The isomorphism classes 
of $\gamma$-indecomposables thus correspond bijectively to $R^+$, and there is 
a partial ordering, again denoted by $\preceq$, on the set of 
$\gamma$-indecomposables defined by
$V_1\preceq V_2$ if there exist indecomposable direct summands $U_1$ of $V_1$
and $U_2$ of $V_2$ such that $U_1\preceq U_2$. That this indeed defines a 
partial ordering (more precisely, that $\preceq$ is anti-symmetric), 
follows easily from the next lemma.

\begin{lemma}\label{strucgind} If $V$ is 
$\gamma$-indecomposable, then ${\rm End}_{k\Gamma}(V)$ is isomorphic to $n$ 
copies of $k$, where $n$ is the symmetrization index of 
$V$, and ${\rm Ext}^1_{k\Gamma}(V,V)=0$.
\end{lemma}

\proof By Lemma \ref{gks}, $V$ is of the form $\widehat{U}$ for
an indecomposable $U\in\bmod k{\Gamma}$. Suppose there exists a 
$k\in\{0,\ldots,n-1\}$ such that $[U,\gamma^kU]\not=0$ or 
$[\gamma^kU,U]^1\not=0$. Then $U\preceq\gamma^kU$ by the properties of the
partial ordering $\preceq$ mentioned above. This leads to a chain
$U\preceq\gamma^kU\preceq\gamma^{2k}U\preceq\ldots\preceq U$. But this implies
$U\simeq\gamma^kU$, contradicting the choice of $n$. Since $U$ fulfills 
$[U,U]=1$ and $[U,U]^1=0$, this proves the lemma. \hb

\begin{lemma}\label{div} If $V$ is $\gamma$-indecomposable
and
$M$
is $\gamma$-sym\-met\-ric, then $[V,V]$ divides $[M,V]$, $[V,M]$, $[M,V]^1$
and $[V,M]^1$. More precisely, if $V=\widehat{U}$, then
$$\frac{[M,\widehat{U}]}{[\widehat{U},\widehat{U}]}=[M,U].$$
\end{lemma}
 
\proof By Lemma \ref{gks}, we have
$V\simeq\widehat{U}=\oplus_{k=0}^{n-1}\gamma^kU$ for some indecomposable
$U\in\bmod k{\Gamma}$ and $n$ as in Definition \ref{symmion}. Thus, we find
$$[M,V]=\sum_{k=0}^{n-1}[M,\gamma^kU]=n\cdot[M,U],$$
since $M$ is $\gamma$-symmetric. By Lemma \ref{strucgind}, we have $[V,V]=n$. 
The other parts of the lemma are proved in the same way. \hb
 
A morphism $f\in{\rm Hom}_{k{\Gamma}}(M,N)={\rm
Hom}_{k\Gamma}(\gamma M,\gamma N)$ between
$\gamma$-symmetric representations $M,N\in\bmod k{\Gamma}$ is called
$\gamma$-symmetric if there exist isomorphisms $\phi:\gamma
M\stackrel{\sim}{\rightarrow} M$
and $\psi:\gamma N\stackrel{\sim}{\rightarrow} N$ such that $f\phi=\psi f$.
It is then easy to
see that kernels, images and cokernels of $\gamma$-symmetric morphisms are
$\gamma$-symmetric.
 
\begin{lemma}\label{exgsm} If $M,N\in\bmod k{\Gamma}$ are $\gamma$-symmetric
such
that $[M,N]\not=0$, then there exists a non-zero
$\gamma$-symmetric morphism from $M$ to $N$.
\end{lemma}

\proof By Lemma \ref{gks}, we can assume without loss of generality that 
$N=\widehat{U}$ for an indecomposable $U$ with symmetrization index $n$. Since 
$[M,N]\not=0$, we also have $[M,U]\not=0$ by Lemma \ref{div}; let 
$f:M\rightarrow U$ be non-zero. Under the isomorphism ${\rm 
Hom}_{k\Gamma}(M,U)={\rm 
Hom}_{k\Gamma}(\gamma^kM,\gamma^kU)\simeq{\rm Hom}_{k\Gamma}(M,\gamma^kU)$, 
the morphism $f$ 
corresponds to a morphism $\gamma^kf:M\rightarrow\gamma^kU$ for each 
$k=0\ldots n-1$. It is then easy to see that the morphism
$$\bigoplus_{k=0}^{n-1}\gamma^kf:M\rightarrow\bigoplus_{k=0}^{n-1}\gamma^kU=N$$
is $\gamma$-symmetric. \hb
 
Given a dimension vector $d=\sum_{i\in\Gamma_0}d_ii\in{\bf N}\Gamma_0$, denote
by $R_d$ the affine space $R_d=\oplus_{\alpha:i\rightarrow j}{\rm
Hom}_k(k^{d_i},k^{d_j})$, and by $G_d$ the group $G_d=\prod_{i\in\Gamma_0}{\rm
GL}(k^{d_i})$. The linear reductive algebraic group $G_d$ acts on the affine
variety $R_d$ by conjugation, i.e.~by
$$(g_i)_i\cdot(M_\alpha)_\alpha=(g_jM_\alpha g_i^{-1})_{\alpha:i\rightarrow
j}.$$
The orbits of $G_d$ in $R_d$ correspond bijectively to the isoclasses of
representations of $\Gamma$ of dimension vector $d$. We denote by ${\cal
O}_M$ the orbit corresponding to the isoclass $[M]$. We say that a
representation $M$ degenerates to $N$ and write $M\leq N$ if the
orbit closure (in the Zariski topology) of the orbit ${\cal O}_M$ contains
the orbit ${\cal O}_N$. This defines a partial ordering on the isoclasses. By
\cite{Bo}, a degeneration $M\leq N$ implies $[U,M]\leq[U,N]$ and $[U,M]^1\leq
[U,N]^1$ for all representations $U$.\\[1ex]
Since there are only finitely many orbits of $G_d$ in the affine space $R_d$,
there has to exist a dense one, whose corresponding representation is denoted
by $E_d$. Thus, we have $E_d\leq
M$ for all representations $M$ of dimension vector $d$. The representation 
$E_d$ is characterized by the property $[E_d,E_d]^1=0$. The next property 
follows immediately.
 
\begin{lemma}\label{symwse} If $d$ is a $\gamma$-symmetric dimension vector,
then the representation $E_d$ is $\gamma$-symmetric.
\end{lemma}
 
\section{The monoid of generic extensions}\label{mon}
 
We continue to use the notation of the previous section. In particular, let
${\Gamma}$ be a quiver of Dynkin type, and let $\gamma$ be an automorphism of
${\Gamma}$. The following
lemmas are proved in \cite{Re1} using the geometry of the representation 
varieties $R_d$.
 
\begin{lemma}\label{exgenext} Given representations $M,N\in\bmod k{\Gamma}$,
there
exists a unique (up to isomorphism) representation $X\in\bmod k{\Gamma}$ such
that
\begin{enumerate}
\item $X$ is an extension of $M$ by $N$, i.e.~there exists an exact sequence
$$\ses{N}{X}{M},$$
\item $\dim_k{\rm End}(X)$ is minimal among all extensions of $M$ by $N$.
\end{enumerate}
\end{lemma}
 
We denote the representation $X$ provided by this lemma by $M*N$, and call it
the generic extension of $M$ by $N$.
 
\begin{lemma} For all representations $L,M,N\in\bmod k{\Gamma}$, we have
$$(L*M)*N\simeq L*(M*N).$$
\end{lemma}
 
Thus, the set of isoclasses $[M]$ of representations $M\in\bmod k{\Gamma}$,
together
with the operation $[M]*[N]=[M*N]$, defines a monoid $\widetilde{\bf \sf M}$
with
unit element $[0]$, the isoclass of the zero representation.
 
\begin{lemma} IF $M,N\in\bmod k{\Gamma}$ are $\gamma$-symmetric, then $M*N$ is
so.
\end{lemma}
 
\proof Applying $\gamma$ to the exact sequence defining $M*N$, we get an exact
sequence
$$\ses{\underbrace{\gamma N}_{\simeq N}}{\gamma(M*N)}{\underbrace{\gamma
M}_{\simeq M}}.$$
Since $\dim_k{\rm End}(\gamma(M*N))=\dim_k{\rm End}(M*N)$, both conditions of
Lemma \ref{exgenext} defining the generic extension of $M$ by $N$ are
fulfilled, thus $\gamma(M*N)\simeq M*N$ by uniqueness. \hb
 
The $\gamma$-symmetric representations of $\bmod k{\Gamma}$ thus form a
submonoid
${\bf \sf M}$ of $\widetilde{\bf \sf M}$.
 
\begin{definition}\label{defmassub} The submonoid ${\bf \sf 
M}({\Gamma},\gamma)={\bf\sf M}\subset\widetilde{\bf\sf M}$ is
called
the
monoid of generic extensions of the pair $({\Gamma},\gamma)$.
\end{definition}
 
\remark The statements of section \ref{state} are proved by showing that
${\bf \sf M}({\Gamma},\gamma)\simeq{\bf \sf U}(C,\leq)$, where $(C,\leq)$ is
the Cartan matrix over a totally ordered index set constructed from 
$({\Gamma},\gamma)$ in Definition
\ref{cartan}.\\[2ex]
In the remaining part of this section, we generalize some results of 
(\cite{Re1}, 3.) to
the monoid ${\bf \sf M}$.\\[2ex]
We enumerate the $\gamma$-indecomposables in $\bmod k{\Gamma}$ as
$V_1,\ldots,V_\nu$,
in such a way that $V_k\preceq V_l$ implies $k\leq l$.
 
\begin{lemma}\label{pbwm} Any element $[M]$ of ${\bf \sf M}$ can be written as
$$[M]=[V_1]^{*m_1}*\ldots*[V_\nu]^{*m_\nu}$$
for certain $m_k\in{\bf N}$ in a unique way.
\end{lemma}
 
\proof By Lemma \ref{gks}, $M$ can be decomposed uniquely as 
$M\simeq\oplus_{k=1}^\nu
V_k^{m_k}$. If $k\leq l$, then $[V_k,V_l]^1=0$, thus $[V_k]*[V_l]=[V_k\oplus
V_l]$. Iterating this, we can write $[M]$ in ${\bf \sf M}$ in the desired
form.
\hb
 
\begin{lemma}\label{gensim} The monoid ${\bf \sf M}$ is generated by the
$\gamma$-symmetrizations
of
the simples in $\bmod k{\Gamma}$, i.e.~by the $[E_{\bf i}]$ for ${\bf i}\in I$.
\end{lemma}
 
\proof Let $M$ be $\gamma$-symmetric. We want to show that $[M]$ can be
written as a product of the $[E_{\bf i}]$ in ${\bf \sf M}$. By the
previous
lemma, we can assume without loss of generality that $M$ is
$\gamma$-indecomposable. In particular, we have $[M,M]^1=0$ by Lemma 
\ref{strucgind}. Let $i$ be a sink
in $\supp M$. Then there exists an
embedding $E_{i}\rightarrow M$, thus a $\gamma$-symmetric embedding
$E_{\bf i}\rightarrow M$ by Lemma \ref{exgsm},
since $M$ is $\gamma$-symmetric. Denoting by $N$ the ($\gamma$-symmetric)
cokernel of this embedding, we have $[M]=[N]*[E_{\bf i}]$, since $M$ has
no
self-extensions. Proceeding by induction on the dimension of $M$ we are done.
\hb
 
\begin{proposition}\label{straight} Let $M,N\in \bmod\,k{\Gamma}$ be
representations
without self-ex\-ten\-sions, and assume that $[N,M]^1=0$. Then $M*N$ has no
self-extensions.
\end{proposition}
 
\proof This is a slight generalization of (\cite{Bo}, Theorem 4.5.). Let
$X$ be the representation $E_{\dimv M+\dimv N}$ without self-extensions of
dimension vector $\dimv
M+\dimv N$. Then $X$ degenerates to $M\oplus N$, and we have $[N,X]^1\leq
[N,M\oplus N]^1=0$ by both (\cite{Bo}, 2.1.) and the assumptions. Thus, we have
$$[N,X]-[N,M\oplus N]=[N,X]^1-[N,M\oplus M]^1=0,$$
using the properties of the Euler form.
Now Theorem 2.4.~of \cite{Bo} shows that we also have an embedding
$N\rightarrow X$, since $N$ embeds into $M\oplus N$, the representation $X$
degenerates to
$M\oplus N$, and $[N,X]=[N,M\oplus N]$. Denoting by $L$ the cokernel of this
embedding, we arrive at the situation
$$\begin{array}{ccccccccc}0&\rightarrow&N&\rightarrow&M*N&\rightarrow&M&\rightarrow&0\\
&&\|&&&&\wedge |&&\\
0&\rightarrow&N&\rightarrow&X&\rightarrow&L&\rightarrow&0.\end{array}$$
By Proposition 2.4.~of \cite{Re1}, this yields a degeneration of $M*N$ to $X$.
Thus $M*N\simeq X$ since $X$ has no self-extensions. \hb
 
\begin{corollary}\label{straightind} Let $M$ and $N$ be 
$\gamma$-indecomposables such that $M\preceq N$.
Then $M*N$ has no self-extensions.
\end{corollary}
 
\proof By Lemma \ref{strucgind}, $\gamma$-indecomposables have no 
self-extensions. By the properties of the partial ordering
$\preceq$,
we have $[N,M]^1=0$. Thus, Proposition \ref{straight} applies. \hb
 
\begin{lemma}\label{eddec} Let ${\bf i}_1,\ldots,{\bf i}_n$ be an enumeration 
of
$I$ such that $k<l$
if there exists an arrow from some vertex in ${\bf i}_k$ to some vertex in 
${\bf i}_l$.
Then for all $d=\sum_{k=1}^nd_k\alpha_{{\bf i}_k}\in Q^+$, we have
$$[E_d]=[E_{{\bf i}_1}]^{*d_1}*\ldots*[E_{{\bf i}_n}]^{*d_n}\mbox{ in }{\bf 
\sf M}.$$
\end{lemma}
 
\proof We just have to note that each representation of dimension vector $d$ 
has a composition series
$$M=M_n\supset\ldots\supset M_1\supset M_0=0$$
which successive subquotients $M_k/M_{k-1}\simeq E_{{\bf i}_k}^{d_k}$ for 
$k=1\ldots n$, and that
$E_d$ has no self-extensions by definition. \hb

\section{Isomorphism of the monoids ${\bf \sf M}$ and ${\bf \sf U}$}\label{iso}
 
In this section, we prove the Realization Theorem \ref{real} in the following
form:
 
\begin{theorem}\label{comp} Let ${\Gamma}$ be a quiver of Dynkin type, and let
$\gamma$ be an
automorphism of ${\Gamma}$. Let $(C,\leq)$ be the pair constructed from
$({\Gamma},\gamma)$
in Definition \ref{cartan}. Then $${\bf \sf U}(C,\leq)\simeq{\bf \sf
M}({\Gamma},\gamma).$$
\end{theorem}

\remark This theorem generalizes Theorem 4.2.~of \cite{Re1}, whose proof 
contains a gap. Namely, it is not clear whether ${\cal U}^+_q$ is isomorphic 
to the generic Hall algebra over ${\bf Q}[[q]]$, as used there.\\[2ex]
We start the proof by constructing a monoid morphism from ${\bf \sf U}={\bf
\sf U}(C,\leq)$ to ${\bf \sf M}={\bf \sf M}({\Gamma},\gamma)$.

\begin{lemma}\label{natdefrel} The defining relations of ${\bf \sf U}$ hold in 
${\bf \sf M}$ if
${\bf i}$ is
replaced by $[E_{\bf i}]$, i.e.
$$[E_{\bf i}]^{*p}*[E_{\bf j}]^{*q}*[E_{\bf i}]^{*r}*[E_{\bf
j}]^{*s}=[E_{\bf i}]^{*(p+r)}*[E_{\bf j}]^{*(q+s)}$$
for ${\bf i}<{\bf j}$ in $I$, and $((p,q),(r,s))$ being two consecutive
entries of the list $L_{\bf i\bf j}$ as
in Definition \ref{defu}.
\end{lemma}
 
\proof Without loss of generality, we can work in a rank $2$ situation, 
i.e.~we can assume
$I=\{{\bf i},{\bf j}\}$. In case $a_{{\bf i}{\bf j}}=0=a_{{\bf j}{\bf i}}$, we
obviously have
$[E_{\bf i},E_{\bf j}]^1=0=[E_{\bf j},E_{\bf i}]^1$, and thus $[E_{\bf
i}]*[E_{\bf j}]=[E_{\bf i}\oplus
E_{\bf j}]=[E_{\bf j}]*[E_{\bf i}]$. So assume that $a_{{\bf i}{\bf
j}}\not=0$. We only treat the cases
where $a_{{\bf i}{\bf j}}=-1$; the other cases can be proved dually and are
left to the
reader. The pair $(C,\leq)$ is then associated to the quiver
$$i_1\rightarrow j_1\mbox{ if }a_{{\bf j}{\bf i}}=-1,\;\;\;
\begin{array}{ccc}i_1&&\\
&\searrow&\\
&&j_1\\
&\nearrow&\\
i_2&&\end{array}\mbox{ if }a_{{\bf j}{\bf i}}=-2,\;\;\;
\begin{array}{ccc}
i_1&&\\
&\searrow&\\
i_2&\rightarrow&j_1\\
&\nearrow&\\
i_3&&\end{array}\mbox{ if }a_{{\bf j}{\bf i}}=-3,$$
respectively, where $\{i_1,\ldots,i_{-a_{{\bf j}{\bf i}}}\}$ forms the
$\gamma$-orbit ${\bf i}$,
and the $\gamma$-orbit ${\bf j}$ consists of the single element $j_1$.\\[1ex]
Calculating the Auslander-Reiten quiver (see \cite{ARS}), we get the following
directed enumerations
of the $\gamma$-indecomposables:
\begin{itemize}
\item $(j_1,i_1j_1,i_1)$ if $a_{{\bf j}{\bf i}}=-1$,
\item $(j_1,i_1j_1\oplus i_2j_1,i_1i_2j_1,i_1\oplus i_2)$ if $a_{{\bf j}{\bf
i}}=-2$,
\item $(j_1,i_1j_1\oplus i_2j_1\oplus i_3j_1,i_1i_2i_3j_1^2,i_1i_2j_1\oplus
i_1i_3j_1\oplus i_2i_3j_1,i_1i_2i_3j_1,i_1\oplus i_2\oplus i_3)$ if
$a_{{\bf j}{\bf i}}=-3$,
\end{itemize}
respectively, where $i_1^{d_1}\ldots i_k^{d_k}j_1^e$ stands for the
indecomposable representation in $\bmod k{\Gamma}$ of dimension vector
$\sum_{l=1}^{k}d_li_l+ej_1$.\\[1ex]
Using the procedure of the proof of Lemma \ref{gensim} or Proposition
\ref{straight}, we see that the elements of ${\bf \sf M}$ corresponding to the
above $\gamma$-indecomposables can be written as
\begin{itemize}
\item $([E_{\bf j}],[E_{\bf i}]*[E_{\bf j}],[E_{\bf i}])$ if
$a_{{\bf j}{\bf i}}=-1$,
\item
$([E_{\bf j}],[E_{\bf i}]*[E_{\bf j}]^{*2},[E_{\bf i}]*[E_{\bf j}],[E_{\bf
i}])$
if $a_{{\bf j}{\bf i}}=-2$,
\item $
([E_{\bf j}],[E_{\bf i}]*[E_{\bf j}]^{*3},[E_{\bf i}]*[E_{\bf j}]^{*2},
[E_{\bf i}]^{*2}*[E_{\bf j}]^{*3},[E_{\bf i}]*[E_{\bf j}],[E_{\bf i}])$
if $a_{{\bf j}{\bf i}}=-3$,
\end{itemize}
respectively. From this, we see that all pairs $({\bf i}^p{\bf j}^q,{\bf
i}^r{\bf j}^s)$ which enter in
the defining relations of ${\bf \sf U}$ correspond in ${\bf \sf M}$ to pairs
of
$\gamma$-indecomposables $(U,V)$ satisfying $U\preceq V$. Thus, Corollary
\ref{straightind} applies and the relations are proved by Lemma \ref{eddec}. 
\hb
 
As a consequence, we get:
 
\begin{corollary}\label{etasur}  The map ${\bf i}\mapsto [E_{\bf i}]$ for
${\bf i}\in
I$ extends to a surjective monoid homomorphism $\eta:{\bf \sf
U}\rightarrow{\bf \sf
M}$.
\end{corollary}
 
\proof The map ${\bf i}\mapsto [E_{\bf i}]$ extends to a monoid homomorphism
since the defining relations of ${\bf\sf U}$ hold in ${\bf \sf M}$. By Lemma
\ref{gensim}, the elements $[E_{\bf i}]$ generate ${\bf\sf M}$. \hb
 
The main difficulty in the proof of Theorem \ref{comp} is to show the
injectivity of the comparison map $\eta$. We first reduce the problem to the
following ``straightening rule", which is the analogue of Proposition 
\ref{straight} in ${\bf\sf U}$:
 
\begin{proposition}\label{straightu} Let $M,N\in \bmod k{\Gamma}$ be
$\gamma$-symmetric representations
without self-ex\-ten\-sions, and assume that $[N,M]^1=0$. Then
$(\dimv M)\cdot(\dimv N)=(\dimv M+\dimv N)$ in ${\bf\sf U}$.
\end{proposition}
 
Using the above comparison map $\eta$, we can now reduce Theorem \ref{comp} to
Proposition \ref{straightu}:
 
\begin{lemma}\label{red1} Theorem \ref{comp} holds provided that Proposition 
\ref{straightu}
holds.
\end{lemma}
 
\proof Recall the enumeration $V_1,\ldots,V_\nu$ of (the isoclasses of) the
$\gamma$-in\-de\-com\-po\-sables from the previous section. Assume that
Proposition \ref{straightu} holds. Let
$$w=(\dimv V_{i_1})\cdot\ldots\cdot(\dimv V_{i_m})$$
be a word in ${\bf\sf U}$. Note that any word in ${\bf\sf U}$ can be written
in this form,
since
$${\bf i}_1\ldots {\bf i}_m=(\dimv E_{{\bf i}_1})\cdot\ldots\cdot(\dimv
E_{{\bf i}_m})$$
using an enumeration ${\bf i}_1,\ldots,{\bf i}_m$ of $I$ as in Lemma 
\ref{eddec}.
We prove that $w$ can be rewritten in the form
$$w=(\dimv V_1)^{n_1}\ldots(\dimv V_\nu)^{n_\nu}.$$
Assume there exists an index $k$ such that $i_k>i_{k+1}$ (otherwise we are
done). Then $[V_{i_{k+1}},V_{i_k}]^1=0$ by the properties of the partial
ordering $\preceq$.
Thus, we can apply Proposition \ref{straightu} to get
$$(\dimv V_{i_k})\cdot(\dimv V_{i_{k+1}})=(\dimv V_{i_k}+\dimv V_{i_{k+1}})$$
in ${\bf\sf U}$. Set $X=E_d$ for $d=\dimv V_{i_k}+\dimv V_{i_{k+1}}$. Writing
$X=\oplus_lV_l^{m_l}$ as in the proof of Lemma \ref{pbwm}, we have
$[V_l,V_{l'}]^1=0$ whenever $m_l\not=0\not=m_{l'}$. Moreover, we have
$i_{k+1}\leq l\leq i_k$ whenever $m_l\not=0$ by $\gamma$-indecomposability of
$V_{i_k}$, $V_{i_{k+1}}$. Applying Proposition
\ref{straightu} again several times, we get $$(\dimv V_{i_k}+\dimv
V_{i_{k+1}})=(\dimv
V_1)^{m_1}\cdot\ldots\cdot(\dimv V_\nu)^{m_\nu}$$
in ${\bf \sf U}$. Putting these two equations together, we arrive at
the following rewriting of $w$:
$$(\dimv V_{i_1})\cdot\ldots\cdot(\dimv
V_{i_{k-1}})\cdot(\dimv V_1)^{m_1}\cdot\ldots\cdot(\dimv V_\nu)^{m_\nu}\cdot 
(\dimv
V_{i_{k+2}})\cdot\ldots\cdot(\dimv V_{i_m}).$$
After a finite number of such rewritings, we obviously arrive at the desired
form
$$w=(\dimv V_1)^{n_1}\cdot\ldots\cdot(\dimv V_\nu)^{n_\nu}.$$
But the image of $w$ under the map $\eta$ is, by definition, the product
$$\eta(w)=[V_1]^{*n_1}*\ldots*[V_\nu]^{*n_\nu},$$
proving the injectivity of $\eta$ by Lemma \ref{pbwm}. We conclude using 
Corollary \ref{etasur}
that the map $\eta$ is a bijection. \hb
 
To prove Proposition \ref{straightu}, we perform a sequence of reductions,
until finally arriving at a rank $2$ situation.
 
\begin{lemma}\label{l56} Proposition \ref{straightu} holds provided it holds
for all
pairs $M,N\in \bmod k{\Gamma}$ of $\gamma$-symmetric representations
without self-ex\-ten\-sions such that $[N,M]^1=0$, but
\begin{itemize}
\item $[N,V]^1\not=0$ for all $\gamma$-symmetric non-zero proper
subrepresentations $V$ of $M$,
\item $[W,M]^1\not=0$ for all $\gamma$-symmetric non-zero proper factor
representations $W$ of $N$.
\end{itemize}
\end{lemma}
 
\proof Let $M,N$ be representations as in Proposition \ref{straightu}. We only
prove the first condition, the second one can be treated dually. We
proceed by induction on the dimension of $M$. Assume
that there exists a $\gamma$-symmetric non-zero proper subrepresentation $V$
of $M$ such that $[N,V]^1=0$, and consider the exact sequence $\ses{V}{M}{X}$.
We apply the functors
${\rm Hom}_{k{\Gamma}}(M,\_)$ and ${\rm Hom}_{k{\Gamma}}(\_,X)$ and get
surjections
$${\rm Ext}_{k{\Gamma}}^1(M,M)\rightarrow {\rm
Ext}_{k{\Gamma}}^1(M,X)\rightarrow {\rm
Ext}_{k{\Gamma}}^1(V,X),$$
hence $[V,X]^1=0$ since $[M,M]^1=0$. Set $M_1=E_{\dimv X}$ and $M_2=E_{\dimv
V}$. Since $V$ is
$\gamma$-symmetric, its dimension vector $\dimv V$ is $\gamma$-symmetric, and
so is $\dimv X=\dimv M-\dimv V$. Hence $M_1$ and $M_2$ are $\gamma$-symmetric
by Lemma \ref{symwse}. We have degenerations $M_1\leq X$ and $M_2\leq V$, thus
$[M_2,M_1]^1\leq[V,X]^1=0$. Thus, we can apply Proposition
\ref{straight} and get $M_1*M_2\simeq M$.\\[1ex]
Denote by $W$ the representation $W=M_2*N$. From the long exact sequence
induced by ${\rm Hom}_{k{\Gamma}}(\_,M_1)$ on the defining exact sequence
$$\ses{N}{W}{M_2},$$
we get $[W,M_1]^1=0$. By induction, we can apply Proposition \ref{straightu}
to the pairs $(M_1,M_2)$, $(M_2,N)$ (since $[N,M_2]^1\leq [N,V]^1=0$ by 
assumption) and
$(M_1,W)$ (since $[W,W]^1=0$ by Proposition \ref{straight}), respectively, and 
get:
$$(\dimv M)\cdot(\dimv N)=(\dimv M_1)\cdot(\dimv M_2)\cdot(\dimv N)=(\dimv
M_1)\cdot(\dimv M_2+\dimv N)$$
$$=(\dimv M_1)\cdot(\dimv W)=(\dimv M_1+\dimv W)=(\dimv M+\dimv N).$$ \hb
 
Note that, in particular, representations $M$ and $N$ as in the lemma have to
be $\gamma$-indecomposable.
 
\begin{lemma}\label{tr} Let $M,N\in \bmod k{\Gamma}$ be $\gamma$-symmetric
representations
without self-ex\-ten\-sions such that $[N,M]^1=0$, and
\begin{itemize}
\item $[N,V]^1\not=0$ for all $\gamma$-symmetric non-zero proper
subrepresentations $V$ of $M$,
\item $[W,M]^1\not=0$ for all $\gamma$-symmetric non-zero proper factor
representations $W$ of $N$.
\end{itemize}
Then at least one of the two following statements holds:
\begin{itemize}
\item $[X,M]^1\not=0$ for all $\gamma$-symmetric non-zero proper factor
representations $X$ of $M$,
\item $[N,Y]^1\not=0$ for all $\gamma$-symmetric non-zero proper
subrepresentations $Y$ of $N$.
\end{itemize}
\end{lemma}
 
\proof Assume that there exists a $\gamma$-symmetric proper factor
representation $X\not=0$ of $M$ such that $[X,M]^1=0$. This defines an exact
sequence $$\ses{V}{M}{X},$$
and we have $[N,V]^1\not=0$ by assumption.
Considering the long exact sequence induced by ${\rm Hom}_{k{\Gamma}}(N,\_)$,
we get
$[N,X]\not=0$ since $[N,M]^1=0$. Let $f:N\rightarrow X$ be a
$\gamma$-symmetric non-zero
morphism provided by Lemma \ref{exgsm}. From the embedding ${\rm Im} f\subset
X$ and the assumption
$[X,M]^1=0$, we deduce $[{\rm Im }f,M]^1=0$. But ${\rm Im}f\not=0$ is a
$\gamma$-symmetric
subrepresentation of $N$, so we have ${\rm Im }f\simeq N$, i.e.~$f$ is
injective.\\[1ex]
If there exists a $\gamma$-symmetric proper subrepresentation $Y\not=0$ of $N$
such that $[N,Y]^1=0$, we can dualize the above argument to construct a
$\gamma$-symmetric non-zero morphism $g:Y\rightarrow M$, which has to be
surjective since
$[N,{\rm Im }g]^1=0$. This would yield a chain of inequalities
$$\dim Y\leq\dim N\leq\dim X\leq\dim M\leq\dim Y,$$
a contradiction. \hb
 
Assume from now on - without loss of generality - that the first case of Lemma
\ref{tr} holds (the other case can be treated dually).
 
\begin{lemma}\label{irr} Assume that $M\in\bmod k{\Gamma}$ is a
$\gamma$-symmetric
representation without self-extensions such that $[X,M]^1\not=0$ for all
$\gamma$-symmetric non-zero proper factor representations $X$ of $M$. Then $M$
is simple, or $\suppb M$ is of type $G_2$, or $M=P_{i_1}$ for the quiver
$$\begin{array}{ccccccc}&&\bullet&\rightarrow&\ldots&\rightarrow&j_1\\
&\nearrow&&&&&\\
i_1&&&&&&\\
&\searrow&&&&&\\
&&\bullet&\rightarrow&\ldots&\rightarrow&j_2,\end{array}$$
and $\gamma$ is of order $2$.
\end{lemma}
 
\proof Assume that $M$ has the above properties; in particular, $M$ is already 
$\gamma$-indecomposable.
Let ${\bf i}$ be an orbit in
the support $\supp M\subset {\Gamma}_0$ of $M$. Then we
can choose a
$\gamma$-symmetric non-zero homomorphism $f:P_{\bf i}\rightarrow M$. If
$f$ is surjective, then ${\bf i}$ has to be the unique source in $\suppb
M\subset I$ by
definition. Otherwise, we get a non-split exact sequence $\ses{{\rm
Im}f}{M}{X}$, such that $X\not=0$ is $\gamma$-symmetric. Thus $[X,M]^1\not=0$
by assumption. Consider the induced exact sequence
$$0\rightarrow{\rm Hom}_{k{\Gamma}}(X,M)\rightarrow{\rm
End}_{k{\Gamma}}(M)\rightarrow
{\rm Hom}_{k{\Gamma}}({\rm Im}f,M)\rightarrow$$
$$\rightarrow{\rm Ext}^1_{k{\Gamma}}(X,M)\rightarrow{\rm
Ext}^1_{k{\Gamma}}(M,M)=0.$$
The image $J$ of ${\rm Hom}_{k{\Gamma}}(X,M)$ in ${\rm End}_{k{\Gamma}}(M)$
consists entirely
of non-invertible $\gamma$-symmetric endomorphisms, since the above short
exact
sequence is non-split. But since the endomorphism ring of each indecomposable
in $\bmod k{\Gamma}$ is trivial, this means that $J=0$, hence $[X,M]=0$. Thus,
we get
$[{\rm Im} f,M]\geq [M,M]+[X,M]^1$. Applying Lemma \ref{div}, this yields an
estimate
$$\frac{[P_{\bf
i},M]}{[M,M]}\geq\frac{[{\rm Im}
f,M]}{[M,M]}\geq\frac{[M,M]}{[M,M]}+\frac{[X,M]^1}{[M,M]}\geq
2.$$
Since $M$ is $\gamma$-indecomposable, we have $M=\widehat{U}$ for some
indecomposable $U\in\bmod k{\Gamma}$ by Lemma \ref{gks}. Again by Lemma 
\ref{div}, we
thus have $\sum_{i\in{\bf i}}\dimv_iU=[P_{\bf i},U]\geq 2$.\\[1ex]
Applying the above argument to all ${\bf i}$ in $\suppb M$, we arrive at one 
of the following two situations:
\begin{enumerate}
\item $\sum_{i\in{\bf i}}\dimv_iU\geq 2$ for all ${\bf i}\in\suppb U$ or
\item There exists a unique source ${\bf i}\in\suppb U\subset I$ such that
$\widehat{U}$
is a factor of $P_{\bf i}$, and $\sum_{i\in {\bf j}}\dimv_iU\geq 2$ for all
${\bf i}\not={\bf j}\in\suppb U$.
\end{enumerate}
We start by analyzing situation a). By a direct inspection of the root systems
of type $A$, $D$ and $E$ (using the classification of possible
automorphisms $\gamma$ of section \ref{quiv}), we conclude that $\dimv U$ has
to be the maximal positive root for the root system
of type $\supp U$,
and the pair $(\supp U,\gamma)$
has to be one of the following:
\begin{itemize}
\item $\supp U$ of type $D_4$, $\gamma$ of order $3$,
\item $\supp U$ of type $E_8$, $\gamma$ trivial,
\item $\supp U$ of type $E_6$, $\gamma$ of order $2$.
\end{itemize}
In particular, $M=\widehat{U}=U$, since the maximal root is always
$\gamma$-symmetric. In the first case, we are done. In the second case,
we choose an immediate successor $X$ of $U$ with respect to the ordering
$\preceq$ on indecomposables in $\bmod k{\Gamma}$. Since $\dimv U$
is the maximal root, $X$ is a proper
factor of $U$, and $[X,U]^1=0$ by the properties of $\preceq$, a
contradiction. In the third case, we consider again the immediate successors
of $U$. Since $U$ is $\gamma$-symmetric, it belongs to the $\tau$-orbit of a 
projective indecomposable $P_i$,
where $i$ is one of the two $\gamma$-fixed vertices of $\supp U$. In any case, 
$U$ has an odd number of
immediate successors, so among them, there is a $\gamma$-symmetric one $X$. 
Argueing as in the
second case, we obtain a contradiction.\\[1ex]
So assume we are in situation b). Note that $\dimv_jP_i$ equals $1$ if there
exists a path from $i$ to $j$ in ${\Gamma}$, and $0$ otherwise. Using this,
we can again proceed by a direct inspection of the root systems, and we arrive
at one of the following situations:
\begin{itemize}
\item $U$ is simple,
\item $\supp U$ of type $D_4$, $\gamma$ of order $3$,
\item $\supp U$ is the quiver
$$\begin{array}{ccccccc}&&\bullet&\rightarrow&\ldots&\rightarrow&\bullet\\
&\nearrow&&&&&\\
i_1&&&&&&\\
&\searrow&&&&&\\
&&\bullet&\rightarrow&\ldots&\rightarrow&\bullet,\end{array}$$
$\gamma$ is of order $2$, and $U=P_{i_1}$.
\end{itemize}
In each of these cases, we are done. \hb

Using this lemma, we can now perform the final reduction.
 
\begin{lemma}\label{red5} Assume that $M$ and $N$ are as in Lemma \ref{tr}. 
Then $\suppb
M\cup \suppb N\subset I$ is at most of rank $2$.
\end{lemma}
 
\proof In the second case of Lemma \ref{irr}, i.e.~$\suppb M$ being of type
$G_2$, there is nothing to prove, since then $I$ has to be of type $G_2$. So
assume we are in the first or the third case. Let ${\bf i}$ be in $\suppb N$.
Then there exists a $\gamma$-symmetric non-zero morphism $P_{\bf i}\rightarrow
N$. If $N$ is not a factor of $P_{\bf i}$, then an argument similar to the one
in the proof of Lemma \ref{irr} shows that $[P_{\bf i},M]\not=0$, which means
${\bf
i}\in\suppb M$. Again as in the proof of Lemma \ref{irr}, we arrive at one of
the following situations:
\begin{itemize}
\item $\suppb N\subset\suppb M$ or
\item there exists a unique source ${\bf i}\in\suppb N$ such that $N$ is a
factor of $P_{\bf i}$, and $\suppb N\setminus\{{\bf i}\}\subset\suppb M$.
\end{itemize}
Dually, we see that
\begin{itemize}
\item $\suppb M\subset\suppb N$ or
\item there exists a unique sink ${\bf i}\in\suppb M$ such that $M$ is a
subrepresentation of $I_{\bf i}$, and $\suppb M\setminus\{{\bf
i}\}\subset\suppb
N$.
\end{itemize}
This analysis gives us enough information to prove the lemma. In case $M$ is
simple, this is obvious. So assume that $M$ is as in the third case of Lemma
\ref{irr}. In particular, $\suppb M$ is of type $C_n$. If $\suppb N$ is not
contained in $\suppb M$, then the second of the above situations applies. Then
${\Gamma}$ has to be one of the following quivers:
$$\begin{array}{ccccc}&&&&3\\ &&&\nearrow&\\ 1&\rightarrow&2&&\\
&&&\searrow&\\ &&&&4\end{array}\;\;\;\mbox{ or
}\;\;\;\begin{array}{ccccccc}&&&&3&\rightarrow&5\\ &&&\nearrow&&&\\
1&\rightarrow&2&&&&\\ &&&\searrow&&&\\ &&&&4&\rightarrow&6,\end{array}$$
$M=P_2$, and $N$ a factor of $P_1$. In both cases, we easily get a
contradiction by direct inspection of the Auslander-Reiten quivers.
If $\suppb N$ is contained in $\suppb M$, then we use the possible situations
for $\suppb M$ to conclude that $\suppb M\setminus\{j_1,j_2\}\subset{\suppb
N}\subset\suppb M$. Again, this makes the possible cases for $M$ and $N$
explicit, and we can conclude that $\suppb M\cup\suppb N$ already is of type
$B_2$ by inspection of the Auslander-Reiten quiver of $\suppb M$. The details 
are left to the reader.  \hb
 
Thus, we only have to study the rank $2$ cases to prove Proposition
\ref{straightu}.\\[2ex]
{\bf Proof of Proposition \ref{straightu}:} The possible rank 2 situations are
listed at the beginning of this section. Using the Auslander-Reiten quivers,
one easily enumerates all possible pairs $(M,N)$ which are as in Lemma 
\ref{tr}. Apart
from trivial relations as $(i)\cdot(j)=(i+j)$, the relations claimed in
Proposition \ref{straightu} are precisely the defining relations of ${\bf \sf
U}$. \hb
 
We conclude that Theorem \ref{comp} is proved.
 
\section{Proofs of the statements of section \ref{state}}\label{proofs}
 
Using the Realization Theorem \ref{comp} we can now easily prove all the
statements of section \ref{state}. We start with the Degeneration Theorem
\ref{dqg}.\\[2ex]
{\bf Proof of Theorem \ref{dqg}:} Denote by ${\cal U}^+_0$ the specialization
${\bf Q}\otimes_{{\bf Q}[[q]]}{\cal U}^+_q$ of ${\cal U}^+_q$ at $q=0$. First
we show that the map $i\mapsto E_i$ extends to an algebra homomorphism
$\theta:{\bf Q}{\bf \sf U}\rightarrow {\cal U}^+_0$. Thus, we have to verify
that the defining relations of ${\bf \sf U}$ hold in ${\cal U}^+_0$.\\[1ex]
Without loss of generality, we can assume to be in the rank $2$ case. In case
$a_{ij}=0=a_{ji}$, there is nothing to prove; so assume by symmetry that
$a_{ij}=-1$.\\[1ex]
In case $a_{ji}=-1$, the $q$-Serre relations of Lemma \ref{qserre} directly
specialize at $q=0$ to the defining relations $iji=i^2j$, $jij=ij^2$ of ${\bf
\sf U}$.\\[1ex]
In case $a_{ji}=-2$, denote by
$$S^+=E_i^2E_j-(q^2+1)E_iE_jE_i+q^2E_jE_i^2,$$
$$S^-=E_iE_j^3-(q^2+q+1)E_jE_iE_j^2+q(q^2+q+1)E_j^2E_iE_j-q^3E_j^3E_i$$
the elements defining the $q$-Serre relations, i.e.~$S^+=0=S^-$ in ${\cal
U}^+_q$. Specializing these elements to $q=0$ gives the relations $iji=i^2j$
and $jij^2=ij^3$. To get the third relation $ij^2ij=i^2j^3$, consider the
element
$$q^{-1}(S^+E_j^2-E_iS^-)=E_iE_jE_iE_j^2+qE_jE_i^2E_j^2-(q^2+q+1)E_iE_j^2E_iE_j+q^2E_iE_j^3E_i.$$
This is a well-defined element of ${\cal U}^+_q$, which by definition equals
zero in this algebra. It specialises to $ijij^2=ij^2ij$ at $q=0$, so we derive
$$ij^2ij=(iji)j^2=(i^2j)j^2=i^2j^3,$$
as desired. Note however that this relation is not a consequence of the first
two relations in ${\bf \sf U}$.\\[1ex]
We proceed similarly in case $a_{ji}=-3$. The calculations get quite involved,
so we only sketch them here. Again, define elements
$$S^+=E_i^2E_j-(q^3+1)E_iE_jE_i+q^3E_jE_i^2,$$
$$S^-=E_iE_j^4-(q^3+q^2+q+1)E_jE_iE_j^3+q(q^4+q^3+2q^2+q+1)E_j^2E_iE_j^2-$$
$$-q^3(q^3+q^2+q+1)E_j^3E_iE_j+q^6E_j^4E_i.$$
Define the following elements inductively; direct calculations show that each
of them in fact belongs to ${\cal U}^+_q$:
$$X=q^{-1}(S^+E_j^3-E_iS^-),$$
$$U=q^{-1}(XE_j-E_iE_jS^-),$$
$$V=q^{-2}E_iX+(q^{-1}+q^{-2})(S^+E_jE_i-S^+E_iE_j-E_iE_jS^+)E_j^2,$$
$$Y=q^{-1}(VE_j+E_iU-2E_i^2E_jS^-+S^+E_iE_j^4).$$
Specializing the elements $S^-,U,X,Y,V,S^+$ at $q=0$ gives (after some tedious
calculations) the six defining relations of ${\bf \sf U}$, in their order of
appearence in the list $L_{ij}$. For example, let us verify the relation 
$E_iE_j^2E_i^2E_j^3=
E_i^3E_j^5$, assuming that the other five
defining relations of ${\bf\sf U}$ hold in ${\cal U}_0^+$. The element $Y$ 
evaluates to
\begin{eqnarray*}
&&(-2q^4-2q^3-3q^2-2q-1)\cdot E_i^2E_j^3E_iE_j^2+(q^2+q)\cdot
E_i^2E_j^2E_iE_j^3+\\
&&(-q^2+q+2)\cdot E_iE_jE_i^2E_j^4
+(2q^5+3q^4+2q^3+2q^2+q)\cdot E_i^2E_j^4E_iE_j+\\
&&(-2q^5-q^4)\cdot E_i^2E_j^5E_i+(q^2-q-1)\cdot E_jE_i^3E_j^4+\\
&&(q+1)\cdot E_jE_i^2E_jE_iE_j^3+(-q-1)\cdot
E_iE_j^2E_i^2E_j^3,\end{eqnarray*}
which proves that $Y$ belongs to ${\cal U}_q^+$. Specializing to $q=0$
yields the relation
\begin{equation}\label{yat0}-E_i^2E_j^3E_iE_j^2+2E_iE_jE_i^2E_j^4-E_jE_i^3E_j^4+
E_jE_i^2E_jE_iE_j^3-E_iE_j^2E_i^2E_j^3=0.\end{equation}
Using the five relations already known, we have
$$E_jE_i^2E_jE_iE_j^3=E_jE_i^2E_iE_j^4=E_jE_i^3E_j^4,$$
\begin{eqnarray*}
E_iE_jE_i^2E_j^4&=&E_iE_jE_iE_iE_j^4=E_i^2E_jE_iE_j^4=E_i^2E_jE_iE_j^3E_j=\\
&=&E_i^2E_iE_j^4E_j=E_i^3E_j^5\end{eqnarray*}
and
$$E_i^2E_j^3E_iE_j^2=E_i^2E_j^3E_iE_jE_j=E_i^3E_j^4E_j=E_i^3E_j^5.$$
Substituting these relations in (\ref{yat0}) gives the desired relation
$$E_iE_j^2E_i^2E_j^3=E_i^3E_j^5.$$
The other relations are treated similarly.\\[2ex]
So we see that $\theta$ extends to an algebra homomorphism. It is obviously
surjective, since ${\cal U}^+_0$ is generated by the elements $E_i$ for $i\in
I$.\\[1ex]
To prove injectivity, we consider the natural $Q^+$-gradings on ${\bf Q}{\bf
\sf U}$, ${\cal U}^+_v$ and ${\cal U}^+_q$, respectively, which are given by
setting the degree of the generating element $i$ (resp.~$E_i$) to $\alpha_i\in
Q^+$. By definition of ${\cal U}^+_q$, we have the
following chain of inequalities for each $d\in Q^+$:
$$\dim_{{\bf Q}(v)}({\cal U}^+_v)_d\leq \dim_{{\bf Q}(q)}{\bf
Q}(q)\otimes_{{\bf Q}[[q]]}({\cal U}^+_q)_d\leq$$
$$\leq\dim_{\bf Q}{\bf Q}\otimes_{{\bf Q}[[q]]}({\cal U}^+_q)_d\leq\dim_{\bf
Q}({\bf Q}{\bf \sf U})_d.$$
Since the quantized enveloping algebra ${\cal U}^+_v$ has a PBW type basis 
(\cite{Lu}, Corollary 40.2.2.),
the leftmost term of the above chain equals the value ${\cal P}(d)$ of
Kostant's partition function at $d$. On the other hand, the rightmost term
equals the number of isoclasses of $\gamma$-symmetric representations of
dimension vector $d$. By Lemma \ref{gks} and the fact that the roots for $C$
correspond to the $\gamma$-symmetrizations of roots for $\widetilde{C}$, the
rightmost term also equals ${\cal P}(d)$. We conclude that equality
holds in each step of the above chain, and that the map $\theta:{\bf Q}{\bf
\sf U}\rightarrow {\cal U}^+_0$ is already an isomorphism. \hb\\[2ex]
{\bf Proof of Lemma \ref{udec}:} The element $d\in Q^+$ corresponds to an 
element
$d\in (\widetilde{Q}^+)^\gamma$. By Lemma \ref{symwse}, the unique
representation without self-extensions $M$ of dimension vector $d$ is
$\gamma$-symmetric, hence has a unique decomposition $M\simeq\oplus_{\alpha\in
R^+}V_\alpha^{a_\alpha}$ into $\gamma$-indecomposables, which correspond to
roots in $R^+$; this yields a decomposition $d=\sum_{\alpha\in
R^+}a_\alpha\alpha$. If $a_\alpha\not=0\not= a_\beta$, then
$$\langle \alpha,\beta\rangle=\langle\dimv V_\alpha,\dimv
V_\beta\rangle=[V_\alpha,V_\beta]-\underbrace{[V_\alpha,V_\beta]^1}_{=0}\geq
0.$$
This proves the existence of the claimed decomposition.\\[1ex]
To prove uniqueness, start with a decomposition $d=\sum_{\alpha\in
R^+}a_\alpha\alpha$ as in the lemma, and define a $\gamma$-symmetric
representation $M=\oplus_{\alpha\in R^+}V_\alpha^{a_\alpha}$. If
$[V_\alpha,V_\beta]^1\not=0$ for some $\alpha,\beta\in R^+$, then
$[V_\alpha,V_\beta]=0$, thus $\langle \alpha,\beta\rangle<0$ by the properties
of the partial ordering $\preceq$. This yields $[M,M]^1=0$. But since a
representation without self-extensions is uniquely determined by its dimension
vector, this proves uniqueness. \hb
 
{\bf Proof of Proposition \ref{str}:} Given $d,e\in Q^+$ as in the 
proposition,
consider the representations without self-extensions $M$ (resp.~$N$) of
dimension vector $d$ (resp.~$e$). Similar to the previous proof, the
assumption that $\langle\beta,\alpha\rangle\geq0$ whenever ${\bf
a}_d(\alpha)\not=0\not={\bf a}_e(\beta)$ translates into the property
$[N,M]^1=0$. Now Proposition \ref{straightu} yields the desired result
$(d)\cdot(e)=(d+e)$ in ${\bf\sf U}$. \hb
 
Finally, we prove the Parametrization Theorem \ref{para}. It was already noted
in \cite{Re4} that the root-theoretic definition \ref{directed} of a directed 
partition is
equivalent to the following representation-theoretic one:\\[1ex]
A partition ${\cal I}_1\cup\ldots\cup{\cal I}_k$ of the set of isoclasses of
$\gamma$-indecomposables is directed if and only if
\begin{enumerate}
\item $[U,V]^1=0$ for all $U,V\in{\cal I}_s$, $1\leq s\leq k$,
\item $[U,V]^1=0=[V,U]$ for all $U\in{\cal I}_s$, $V\in{\cal I}_t$, $1\leq
s<t\leq k$.
\end{enumerate}
 
Now Lemma \ref{ede} follows directly from the existence of the partial
ordering $\preceq$ on $\gamma$-indecomposables.\\[2ex]
{\bf Proof of Theorem \ref{para}:} Given a $\gamma$-symmetric representation 
$M$
of ${\Gamma}$, we decompose it as $M=M_1\oplus\ldots\oplus M_k$, where $M_s$
denotes the direct sum of all $\gamma$-indecomposable direct summands
$U\in{\cal I}_s$  of $M$. By the definition of a directed partition, we have
$[M_t,M_s]^1=0$ for all $1\leq s\leq t\leq k$. This yields the identity
$[M]=[M_1]*\ldots*[M_k]$ in ${\bf \sf M}$ by definition. Using Lemma
\ref{eddec}, the right hand side of this identity translates under the
isomorphism ${\bf \sf U}\simeq{\bf \sf M}$ into the element $(\dimv
M_1)\cdot\ldots\cdot(\dimv M_k)$. Now these elements constitute precisely the
image of the map $p_{{\cal I}_*}:{\bf N}R^+\rightarrow {\bf\sf U}$ of the 
theorem, whereas
the elements $[M]$ constitute precisely the elements of ${\bf\sf
M}\simeq{\bf\sf U}$. We see that the theorem is proved. \hb

\end{document}